\documentclass[12pt,oneside,reqno]{amsart}
\usepackage[utf8]{inputenc}
\usepackage{amsmath,mathrsfs,amssymb,multirow,setspace,enumerate}
\usepackage[a4paper, total={7in, 9in}]{geometry}
\usepackage{xcolor}
\usepackage{lettrine}
\usepackage{tikz}
\usepackage[11pt]{moresize}
\theoremstyle{plain}
\usepackage{graphicx}
\usepackage{amsmath}
\newtheorem{theorem}{Theorem}[section]
\newtheorem{lemma}[theorem]{Lemma}
\newtheorem{proposition}[theorem]{Proposition}
\newtheorem{corollary}[theorem]{Corollary}

\theoremstyle{definition}

\theoremstyle{remark}
\newtheorem{remark}[theorem]{Remark}

\input xy
\xyoption{all}
\def\c{\mathcal{P}_e(SD_{8n})}
\def\d{\mathcal{P}_e(Q_{4n})}
\def\a{\mathcal{P}_e(D_{2n})}

\begin{document}
\title[Enhanced Power Graph of Certain Non-abelian Groups]{Enhanced Power Graph of Certain Non-abelian Groups}
\author[Jitender Kumar,Sandeep Dalal, and Parveen ]{ Parveen, Sandeep Dalal, and Jitender Kumar }
\address{Department of Mathematics, Birla Institute of Technology  and Science Pilani, Pilani, India}
\email{p.parveenkumar144@gmail.com,deepdalal10@gmail.com,jitenderarora09@gmail.com}
\begin{abstract}
 The enhanced power graph of a group $G$ is a simple undirected graph with vertex set $G$ and two vertices are adjacent if they belong to same cyclic subgroup. In this paper, we study distant properties and detour distant properties such as closure, interior, distance degree sequence and eccentric subgraph of the enhanced power graph of semidihedral group. Consequently, we obtained the metric dimension and resolving polynomial of the enhanced power graph of semidihedral group. At the final part of this paper, we obtained the Laplacian spectrum of the enhanced power graph of semidihedral, dihedral and generalized quaternion groups.
\end{abstract}

\subjclass[]{05C25, 05C50}

\keywords{Enhanced power graph, Resolving polynomial, Laplacian spectrum, Semidihedral group}

\maketitle

 \section{Introduction}
There are many graphs associated with group structure viz. power graphs \cite{a.MKsen2009}, enhanced power graph \cite{a.Bera2017}, commuting graphs \cite{a.segev2001commuting}, cyclic graphs \cite{abdollahi2009noncyclic}, divisibility graphs \cite{a.Adeleh2017}, Cayley graph \cite{kelarev2002undirected,luo2021perfect} etc. The investigation of graphs related to various algebraic structures is very important, because graphs of this type have valuable applications and are related to automata theory (see \cite{kelarev2009cayley,kelarev2004labelled} and the books \cite{kelarev2003graph,kelarev2002ring}). A significant number of papers devoted to algebraic graphs over certain finite groups such as power graphs over dihedral groups, finite cyclic groups, semidihedral groups, dicyclic groups (see \cite{a.2017ashrafi_automorphism,a.chattopadhyay2018connectivity,a.chattopadhyayspectralradius2017,a.chattopadhyay-laplacian,a.2017Hamzaeh_automorphism}); reduced power graphs over dihedral group, semidihedral group and dicyclic group \cite{a.Rajkumar-Laplacian};  commuting graphs over semidihedral group \cite{kumar2021commuting}, generalized dihedral group \cite{a.Kakkar2018commuting}, dicyclic group \cite{chen2020commuting}, dihedral group \cite{a.Ali-2016}; co-prime order graph of dihedral group and dicyclic group \cite{ma2020co-prime}; Cayley graphs over semidihedral group \cite{wang2021prettyCayley}, dihedral group \cite{cao2021perfect} etc. 

 In 2000, Kelarev and Quinn defined two interesting classes of directed graphs, viz.  divisibility and power graph on semigroups \cite{a.kelarev2002directed,a.kelarev2001powermatrices}. The undirected power graph $\mathcal{P}(S)$ of a semigroup $S$  became the main focus of study in \cite{a.MKsen2009} defined by Chakraborty et al. is whose vertex set is $S$ and two distinct vertices $x, y$ are adjacent if either $x = y^m$ or $y = x^n$  for some $m, n \in \mathbb N$. We refer the reader for more results on power graph to survey paper \cite{kumar2021recent}.  The commuting graph $\Delta(S)$ of a semigroup $S$ is the graph whose vertex set is $\Omega \subseteq S$ and two distinct vertices $x, y$ are adjacent if $xy = yx$. The power graph $\mathcal P(S)$ is a spanning subgraph of commuting graph $\Delta (S)$ when $\Omega = S$. In \cite{a.Aalipour2017}, Aalipour et al. characterized the finite group $G$ such that the $\mathcal P(G)$ coincides with the commuting graph of $G$. If these two graphs of $G$ do not coincide, then to measure how much the power graph is close to the commuting graph of a group $G$, they introduced a new graph called \emph{enhanced power graph,} denoted by $\mathcal{P}_e(G)$, is the graph whose vertex set is the group $G$ and two distinct vertices $x, y$ are adjacent if $x, y \in \langle z \rangle$ for some $z \in G$. Aalipour et al. \cite{a.Aalipour2017} characterize the finite group $G$, for which equality holds for either two of the three graphs viz. power graph, enhanced power graph and commuting graph of $G$. Further, the enhanced power graphs have been studied by various researchers. In \cite{a.Bera2017}, Bera et al. characterized the abelian groups and the non abelian $p$-groups having dominatable enhanced power graphs. In \cite{a.Dupont2017}, Dupont et al. determined the rainbow connection number of enhanced power graph of a finite group $G$. Later, Dupont et al. studied the graph theoretic properties in \cite{a.Dupont2017quotient} of enhanced quotient graph of a finite group $G$.  Ma et al. \cite{2019Mametric} investigated the metric dimension of an enhanced power graph of finite groups.  Zahirovi$\acute{c}$ et al. \cite{a.2019study}  proved that two finite abelian groups are isomorphic if their enhanced power graphs are isomorphic. Also, they supplied a characterization of finite nilpotent groups whose enhanced power graphs are perfect.  Recently, Panda et al. \cite{a.Panda-enhanced} studied the graph-theoretic properties viz. minimum degree, independence number, matching number, strong metric dimension and perfectness of enhanced power graph over finite abelian groups and some non abelian groups such as Dihedral groups, Dicyclic groups and the group $U_{6n}$.  Dalal et al. \cite{a.dalal2021enhanced} investigated the graph-theoretic properties of enhanced power graphs over dicyclic group and the group $V_{8n}$.  Bera et al. \cite{a.Bera2021EPG} gave an upper bound for the vertex connectivity of enhanced power graph of any finite abelian group $G$. Moreover, they classified the finite abelian group $G$ such that their proper enhanced power graph is connected.

Spectra of graphs have many extensive implementations in quantum chemistry, though resolving set plays a significant role in the modelling of virus propagation in the computer networks. Moreover, distance degree sequence is a remarkable tool to keep the concealment of exclusive data in some databases. In this paper, we study the enhanced power graph of certain non-abelian groups. In this connection, we have investigated distant properties and detour distant properties of the enhanced power graph of semidihedral group. Moreover, we obtained the Laplacian spectrum of the enhanced power graph of the semidihedral group, generalized quaternion group and dihedral group, respectively. This paper is structured as follows.  In Section $2$, we provide necessary background material and fix our notations used throughout the paper. Section $3$ comprises the results of distant properties and detour distant properties including metric dimension, resolving polynomial etc. The Laplacian spectrum of the enhanced power graph of semidihedral, generalized quaternion and dihedral groups is obtained in Section $4$. 
%In Section $3$, we study some properties of $\Delta(G)$. Section $4$ comprises the study of various graph invariants of $\c$ viz. Hamiltonian, perfectness, independence number, clique number, vertex connectivity, edge connectivity, vertex covering number, edge covering number etc.. Moreover, we study the Laplacian spectrum, metric dimension, resolving polynomial and the detour properties  of the commuting graph of $SD_{8n}$. 

\section{Preliminaries}
In this section, we recall  necessary definitions, results and notations of graph theory from \cite{b.West}.
A graph $\Gamma$ is a pair  $ \Gamma = (V, E)$, where $V = V(\Gamma)$ and $E = E(\Gamma)$ are the set of vertices and edges of $\Gamma$, respectively. Moreover, the \emph{order} of a graph $\Gamma$ is the number of vertices in $\Gamma$. We say that two different vertices $a, b$ are $\mathit{adjacent}$, denoted by $a \sim b$, if there is an edge between $a$ and $b$. Also, we denote this edge $a \sim b$ by $(a, b)$.  The \emph{neighbourhood} $ N(x) $ of a vertex $x$ is the set all vertices adjacent to $x$ in $ \Gamma $. Additionally, we denote $N[x] = N(x) \cup \{x\}$. It is clear that we are considering simple graphs, i.e. undirected graphs with no loops  or repeated edges. A \emph{subgraph} of a graph $\Gamma$ is a graph $\Gamma'$ such that $V(\Gamma') \subseteq V(\Gamma)$ and $E(\Gamma') \subseteq E(\Gamma)$. A graph $\Gamma$ is said to be $complete$ if any two distinct vertices are adjacent. We denote by $K_n$ as the complete graph on $n$ vertices. The $complement \ \overline{\Gamma}$ of a simple graph $\Gamma$ is a simple graph with vertex set $V(\Gamma)$ defined by $(u,v)\in E(\overline{\Gamma})$ if and only if $(u,v)\notin E(\Gamma)$.   A graph with no cycle is called $acyclic$. A $tree$ is a connected acyclic graph. A $spanning \ subgraph$  of $\Gamma$ is a subgraph with vertex set $V(\Gamma)$. A $spanning \ tree$ is a spanning subgraph that is a tree.
A \emph{walk} $\lambda$ in $\Gamma$ from the vertex $u$ to the vertex $w$ is a sequence of  vertices $u = v_1, v_2,\ldots, v_{m} = w$ $(m > 1)$ such that $v_i \sim v_{i + 1}$ for every $i \in \{1, 2, \ldots, m-1\}$. If no edge is repeated in $\lambda$, then it is called a \emph{trail} in $\Gamma$. A trail whose initial and end vertices are identical is called a \emph{closed trail}. A walk is said to be a \emph{path} if no vertex is repeated. The length of a path is the number of edges it contains.

Let $\Gamma$ be a graph. The (\emph{detour}) \emph{distance}, $(d_D(u,v))$ $d(u, v)$, between two vertices $u$ and $v$ in $\Gamma$ is the length of (longest) shortest $u-v$ path in $\Gamma$. The (\emph{detour}) \emph{eccentricity} of a vertex $u$, denoted by $(ecc_D(u))$ $ecc(u)$, is the maximum (detour) distance between $u$ and any vertex of $\Gamma$. The minimum (detour) eccentricity among the vertices of $\Gamma$ is called the (\emph{detour}) \emph{radius} of $\Gamma$, it is denoted by $(rad_D(\Gamma))$ $rad(\Gamma)$. The \emph{detour diameter} of a graph $\Gamma$ is the maximum detour eccentricity in $\Gamma$, denoted by $diam_D(G)$. A vertex $v$ is said to be \emph{eccentric vertex} for $u$ if $d(u, v) = ecc(u)$. A vertex $v$ is said to be an eccentric vertex of $\Gamma$ if $v$ is an eccentric vertex for some vertex of $\Gamma$. A graph $\Gamma$ is said to be an \emph{eccentric graph} if every vertex of $\Gamma$ is an eccentric vertex. The \emph{centre} of $\Gamma$ is a subgraph of $\Gamma$ induced by the vertices having minimum eccentricity and it is denoted by $Cen(\Gamma)$. The \emph{closure} of $\Gamma$ of order $n$ is the graph obtained from $\Gamma$ by recursively joining pairs of non-adjacent vertices whose sum of degree is at least $n$ until no such pair remains and it is denoted by $Cl(\Gamma)$.  The graph $\Gamma$ is said to be \emph{closed} if $\Gamma = Cl(\Gamma)$ (see \cite{b.chartrand2004introduction}).

A vertex $v$ in $\Gamma$ is a \emph{boundary vertex} of a vertex $u$ if $d(u, w) \leq d(u, v)$ for $w \in$ N$(v)$, while a vertex $v$ is a boundary vertex of a graph $\Gamma$ if $v$ is a boundary vertex of some vertex of  $\Gamma$. The subgraph $\Gamma$ induced by its boundary vertices is the \emph{boundary} $\partial(\Gamma)$ of $\Gamma$. A vertex $v$ is said to be a \emph{complete vertex} if the subgraph induced by the neighbors of $v$ is complete. A vertex $v$ is said to be an \emph{interior vertex} of $\Gamma$ if for each $u \ne v$, there exists a vertex $w$  and a path $u-w$ such that $v$ lies in that path at the same distance from both $u$ and $w$. A subgraph induced by the interior vertices of $\Gamma$ is called \emph{interior} of $\Gamma$ and it is denoted by $Int(\Gamma)$. 

\begin{theorem}[\cite{b.chartrand2004introduction}, p.337]\label{SD_8n-boundary-complete}
Let $\Gamma$ be a connected graph and $v \in V(\Gamma)$. Then $v$ is a complete vertex of $\Gamma$  if and only if $v$ is a  boundary vertex of $x$ for all $x \in V(\Gamma) \setminus \{v\}$.
\end{theorem}

\begin{theorem}[\cite{b.chartrand2004introduction}, p.339]\label{SD_8n-boundary-interior}
Let $\Gamma$ be a connected graph and $v \in V(\Gamma)$. Then $v$ is a  boundary vertex of $\Gamma$  if and only if $v$ is not an interior vertex of $\Gamma$.
\end{theorem}

For a finite simple undirected graph $\Gamma$ with vertex set $V(\Gamma) = \{v_1, v_2, \ldots, v_n\}$, the \emph{adjacency matrix} $A(\Gamma)$ is the $n\times n$ matrix with $(i, j)th$ entry is $1$ if $v_i$ and $v_j$ are adjacent and $0$ otherwise. We denote the diagonal matrix $D(\Gamma) = {\rm diag}(d_1, d_2, \ldots, d_n)$ where $d_i$ is the degree of the vertex $v_i$ of $\Gamma$, $i = 1, 2, \ldots, n$. The \emph{Laplacian matrix} $L(\Gamma)$ of $\Gamma$ is the matrix $D(\Gamma) - A(\Gamma)$. The matrix $L(\Gamma)$ is symmetric and positive semidefinite, so that its eigenvalues are real and non-negative. Furthermore, the sum of each row (column) of $L(\Gamma)$ is zero. Recall that, the \emph{characteristic polynomial} of $L(\Gamma)$ is denoted by $\Phi(L(\Gamma), x)$. The eigenvalues of $L(\Gamma)$ are called the \emph{Laplacian eigenvalues} of $\Gamma$ and it is denoted by $\lambda_1(\Gamma) \geq \lambda_2(\Gamma) \geq \cdots \geq \lambda_n(\Gamma) = 0$. Now let $\lambda_{n_1}(\Gamma) \geq \lambda_{n_2}(\Gamma) \geq \cdots \geq \lambda_{n_r}(\Gamma) = 0$ be the distinct eigenvalues of $\Gamma$ with multiplicities $m_1, m_2, \ldots, m_r$, respectively. The \emph{Laplacian spectrum} of $\Gamma$, that is, the spectrum  of $L(\Gamma)$ is represented as $\displaystyle \begin{pmatrix}
\lambda_{n_1}(\Gamma) & \lambda_{n_2}(\Gamma) & \cdots& \lambda_{n_r}(\Gamma)\\
 m_1 & m_2 & \cdots & m_r\\
\end{pmatrix}$. 
We denote the matrix $J_n$ as the square matrix of order $n$ having all the entries as $1$ and $I_n$ is the identity matrix of order $n$.

\section{Distant properties and detour distant properties} 
In this section, we study the distant properties and detour distant properties such as closure, interior, distance degree sequence and eccentric subgraph of the enhanced power graph of semidihedral group. Also, we obtain the metric dimension and resolving polynomial of the enhanced power graph of semidihedral group. For $n \geq 2$, the \emph{semidihedral group} $SD_{8n}$ is a group of order $8n$ is defined in terms of generators and
relations as $$SD_{8n} = \langle a, b  :  a^{4n} = b^2 = e,  ba = a^{2n -1}b \rangle.$$
 We have $$ ba^i = \left\{ \begin{array}{ll}
a^{4n -i}b & \mbox{if $i$ is even,}\\
a^{2n - i}b& \mbox{if $i$ is odd,}\end{array} \right.$$
so that every element of $SD_{8n} {\setminus} \langle a \rangle$ is of the form $a^ib$ for some $0 \leq i \leq 4n-1$. We denote the subgroups $P_i = \langle a^{2i}b \rangle = \{e, a^{2i}b\}$ and $ Q_j =  \langle a^{2j + 1}b \rangle = \{e, a^{2n}, a^{2j +1}b, a^{2n + 2j +1}b\} $. Then we have
\begin{equation}\label{SD_8n}
    SD_{8n} = \langle a \rangle \cup \left( \bigcup\limits_{ i=0}^{2n-1} P_i \right) \cup \left( \bigcup\limits_{ j=  0}^{n-1} Q_{j}\right).
\end{equation}

A \emph{maximal cyclic subgroup} of a group $G$ is a cyclic subgroup of $G$ that is not properly contained in any other cyclic subgroup of $G$.
In the following lemma, we obtain the neighbourhood of all the vertices of $\mathcal{P}_e(SD_{8n})$.

\begin{lemma}\label{nbd}
In  $\mathcal{P}_e(SD_{8n})$, we have
\begin{enumerate}[\rm (i)]
\item {\rm N}$[e] = SD_{8n}$.
\item {\rm N}$[a^{2n}] = \langle a \rangle \cup \{a^{2i + 1}b : 0 \leq i \leq 2n - 1  \}$.

\item {\rm N}$[a^{i}] = \langle a \rangle$, where $1 \leq i \leq 4n - 1$ and $i \ne 2n$.

\item {\rm N}$[a^{2i + 1}b] = \langle a^{2i + 1}b \rangle = \{e, a^{2n}, a^{2i +1}b, a^{2n + 2i +1}b\} $, where $0 \leq i \leq 2n - 1$.

\item {\rm N}$[a^{2i}b] = \{e,a^{2i}b \} $, where $1 \leq i \leq 2n$.
\end{enumerate}
\end{lemma}

\begin{proof}
The proof of (i) is straightforward.\\
(ii) Note that $a^{2n} \sim x$ for all $x \in \langle a \rangle$ and $a^{2n} \sim a^{2i + 1}b$ for all $i$, where $0 \leq i \leq 2n -1$ as $a^{2n} \in \langle a^{2i + 1}b \rangle$. This implies that $\langle a \rangle \cup \{a^{2i + 1}b : 0 \leq i \leq 2n - 1  \} \subseteq {\rm N}[a^{2n}]$. If possible, let $x \in {\rm N}[a^{2n}]$ such that $x = a^{2i}b$ for some $i$, where $0 \leq i \leq 2n - 1$. In view of Equation (\ref{SD_8n}), observe that $x$ belongs to exactly one cyclic subgroup $P_i$ and $x \sim a^{2n}$ gives $a^{2n} \in P_i$; a contradiction.\\
(iii) Let $i \ne 2n$ and $1 \leq i \leq 4n - 1$. Then clearly $\langle a \rangle \subseteq {\rm N}[a^i]$. If $a^i \sim x$ for some $x \in SD_{8n} \setminus \langle a \rangle$, then either $x, a^i \in P_j$ or $x, a^i \in Q_k$ for some $j, k$, where $0 \leq j \leq 2n -1$ and $0 \leq k \leq n - 1$ which is not possible.
Further, note that if $\langle x \rangle$ is a maximal cyclic subgroup, then ${\rm N}[x] = \langle x \rangle$. Since $P_i$ and $Q_j$ are maximal cyclic subgroup generated by $a^{2i}b$ and $a^{2j + 1}b$, respectively,  where $0 \leq i \leq 2n -1$ and $0 \leq j \leq n - 1$. It follows that  {\rm N}$[a^{2i}b] = \{e,a^{2i}b \} $  and ${\rm N}[a^{2j + 1}b] = \langle a^{2j + 1}b \rangle = \{e, a^{2n}, a^{2j +1}b, a^{2n + 2j +1}b\}$. Thus, the result holds.
\end{proof}

Now we obtain the detour radius, detour eccentricity, detour degree, detour degree sequence and detour distance degree sequence of each vertex of $\c$.

\begin{theorem} \label{detour-ecentricity}
In  $\mathcal{P}_e(SD_{8n})$, we have

 \[ecc_D(v) =\left\{ \begin{array}{ll}
4n + 1, & \; \; \; \mbox{when $v \in \{e, a^{2n}\}$;}\\
4n + 2, & \; \; \; \mbox{when $v \in  \{a^{2i}b : 1 \leq i \leq 2n\} $;}\\
4n + 3, & \; \; \; \mbox{when $v \in \{a^i : 1 \leq i \leq 4n - 1, i \ne 2n\} \cup\{a^{2i + 1}b : 1 \leq i \leq 2n - 1  \} $.}
\end{array}\right.\]
\end{theorem}

\begin{proof}
 In view of Lemma \ref{nbd}, we have the following:
 \begin{itemize}
\item $d_D(e,  a^i) = \left\{ \begin{array}{ll}
4n + 1, & \; \; \; \text{if} ~ i \ne 2n\\
4n - 1, & \; \; \; \text{if} ~ i = 2n.
\end{array}\right.$ 

\item $d_D(e, a^{2i + 1}b) = 4n + 1$ for all $i$, where $0 \leq i \leq 2n -1$;

\item $d_D(e, a^{2i}b) =1$ for all $i$, where $1 \leq i \leq 2n$;

\item $d_D(a^{2n}, a^j) = 4n + 1$ for all $j$, where $1 \leq j \leq 4n -1$ and $j \ne 2n$; 

\item $d_D(a^{2n}, a^{2i + 1}b) = 4n + 1$ for all $i$, where $0 \leq i \leq 2n -1$;

\item $d_D(a^{2n}, a^{2i}b) =4n$ for all $i$, where $1 \leq i \leq 2n$;

\item $d_D(a^i, a^j) = 4n +1$ for all $i \ne j$, where $1 \leq i, j \leq 4n  - 1$ and $i, j \ne 2n$;

\item $d_D(a^i, a^{2j + 1}b) = 4n +3$, where $1 \leq i\leq 4n  - 1$, $i\ne 2n$ and $0 \leq j \leq 2n - 1$;

\item $d_D(a^i, a^{2j}b) = 4n +2$, where $1 \leq i\leq 4n  - 1$, $i\ne 2n$ and $0 \leq j \leq 2n$;

\item $d_D(a^{2i + 1}b, a^{2n + 2i + 1}b) = 4n +1$, where $0 \leq i\leq n - 1$;

\item $d_D(a^{2i + 1}b, a^{2j + 1}b) = 4n +3$ for all $i\ne j$, where $0 \leq i\leq 2n  - 1$ and $j \ne n + i$;

\item $d_D(a^{2i + 1}b, a^{2j}b) = 4n +2$, where $0 \leq i\leq 2n  - 1$ and $0 \leq j \leq 2n - 1$;

\item $d_D(a^{2i}b, a^{2j}b) = 2$ for all $i \neq j$, where $1 \leq i, j\leq 2n$.
\end{itemize}
Thus we conclude the result. 
\end{proof}

By the definition of rad$_D(\c)$ and diam$_D(\c)$, we have the following corollary.

\begin{corollary}\label{Detour-radius-diameter}
In $\c$, we have
\begin{enumerate}[\rm(i)]

\item {\rm rad}$_D(\mathcal{P}_e(SD_{8n}) = 4n + 1$

\item {\rm diam}$_D(\mathcal{P}_e(SD_{8n}) = 4n + 3$.
\end{enumerate}
\end{corollary}

The \emph{detour degree} $d_D(v)$ of $v$ is the number $|D(v)|$, where $D(v) = \{u \in V(\Gamma) : d_D(u, v) = ecc_D(v) \}$. The \emph{average detour degree} $D_{av}(\Gamma)$ of a
graph $\Gamma$ is the quotient of the sum of the detour degrees of all the vertices
of $\Gamma$ and the order of $\Gamma$. The detour degrees of the vertices of a graph
written in non-increasing order is said to be the \emph{detour degree sequence} of graph
$\Gamma$, denoted by $D(\Gamma)$. For a vertex $x \in V(\Gamma)$, we denote $D_i(x)$ be the number of vertices at a detour distance $i$ from the vertex $x$. Then the sequence $D_0(x), D_1(x), D_2(x), \ldots, D_{ecc_D(x)}(x)$ is called \emph{detour distance degree sequence of a vertex x}, denoted by $dds_D(x)$. In the remaining part of this paper, $(a^r, b^s, c^t)$ denote $a$ occur $r$ times, $b$ occur $s$ times and $c$ occur $t$ times in the sequence. Now we have the following remark.

\begin{remark}[{\cite[Remark 2.6]{a.Ali-2016}}] In a graph $\Gamma$, we have
\begin{enumerate}[\rm (i)]
\item $D_0(v) = 1$ and $D_{ecc_D}(v) = d_D(v)$.
\item The length of sequence $dds_D(v)$ is one more than the detour eccentricity of $v$.
\item $\displaystyle \sum\limits_{i = 0}^{ecc_D(v)} D_i(v) = |\Gamma|$.
\end{enumerate}
\end{remark}

\begin{proposition}\label{Detour-degree}
In $\c$, we have 

\[d_D(x) =  \left\{ \begin{array}{ll}
6n-4, & \; \; \; \mbox{when $x = a^{2i + 1}b$, where $0 \leq i \leq 2n - 1$};\\
2n, & \; \; \; \mbox{when $x = a^i$ for all $i \ne 2n$, where $1 \leq i \leq 4n - 1$};\\
6n - 2, & \; \; \; \mbox{when $x \in \{ a^{2i}b : 0 \leq i \leq 2n - 1\} \cup \{e,  a^{2n}\}$}.
\end{array}\right.\]
\end{proposition}

\begin{proof} 
Let  $x = a^{2i + 1}b$ for some $i$, where $0 \leq i \leq 2n - 1$.  In view of Theorem \ref{detour-ecentricity}, $ecc_D(x) = 4n +3$. By the proof of Theorem \ref{detour-ecentricity}, one can observe that  \[D(x) = \{a^{2j +1}b : 0 \leq j \leq 2n - 1, j \ne i, j \ne n + i \}\cup \left( \langle a \rangle \setminus \{e, a^{2n}\} \right).\]
For $x = a^i$, where $1 \leq i \leq 4n -1$ and $i \ne 2n$, we have $ecc_D(x) = 4n +3$ (see Theorem \ref{detour-ecentricity}). Again by the proof of Theorem \ref{detour-ecentricity}, we get $D(x) = \{a^{2j +1}b : 0 \leq j \leq 2n - 1\}$.
Now let $x \in \{ a^{2i}b : 0 \leq i \leq 2n - 1\} \cup \{e, a^{2n}\}$. Then $ecc_D(x) = 4n + 1$ when $x \in \{e, a^{2n}\}$, and $ecc_D(x) = 4n + 2$ when $x = a^{2i}b$ for some $i$, where $0 \leq i \leq 2n - 1$ (cf. Theorem \ref{detour-ecentricity}). By the proof of Theorem \ref{detour-ecentricity}, we have \[D(x) = \{a^{2i + 1}b : 0 \leq i \leq 2n - 1\} \cup \left( \langle a \rangle \setminus \{e, a^{2n}\} \right).\]
Thus the result holds.
%Similar to $x \in Z(SD_{8n})$, for $x \in \langle a \rangle \setminus Z(SD_{8n})$ we obtain $D(x) = SD_{8n} \setminus \langle a \rangle$ (cf. Theorem \ref{detour-ecentricity}). Now let $x = a^ib$ for some $i$, where $1 \leq i \leq 4n$. Similar to $x \in Z(SD_{8n})$, when $n$ is even, we get
%\[D(a^ib) = \left(\{a^jb : 1\leq j \leq 4n\} \setminus \{a^ib, a^{2n +i}b\}\right) \cup \left(\langle a\rangle \setminus Z(SD_{8n}) \right).\]
%and for odd $n$,
%\[D(a^ib) = \left(\{a^jb : 1\leq j \leq 4n\} \setminus \{a^{n+i}b, a^{2n +i}b, a^{3n +i}b, a^{4n + i}b\}\right) \cup \left(\langle a\rangle \setminus Z(SD_{8n}) \right).\]
%  \qed
\end{proof}

\begin{corollary} In $\c$, we have

\begin{enumerate} [\rm (i)]
\item \[D(\c) = ((6n-4)^{2n}, (2n)^{4n- 2}, (6n - 2)^{2n + 2}) \]

\item \[D_{av}(\c) = \frac{8n^2-n- 1}{2n} \]
\end{enumerate}
\end{corollary}

\begin{theorem}
In $\c$, we have
\[dds_D(\c) = \left\{ \begin{array}{ll}
(1, 2n, 0^{4n - 3}, 1, 0, 6n -2), (1, 0^{4n -2},1, 2n, 6n -2), (1, 0^{4n}, 4n - 1, (2n)^2)^{4n - 2},\\ 
\vspace{0.1cm}\\
(1^2, 2n - 1, 0^{4n - 3}, 1, 0, 6n -2)^{2n}, (1, 0^{4n}, 3,2n, 6n - 4)^{2n}.
\end{array}\right. \]
\end{theorem}

\begin{proof}
In view of Theorem \ref{detour-ecentricity}, we get
 \[ecc_D(v) =\left\{ \begin{array}{ll}
4n + 1, & \; \; \; \mbox{when $v \in \{e, a^{2n}\}$;}\\
4n + 2, & \; \; \; \mbox{when $v \in  \{a^{2i}b : 1 \leq i \leq 2n\} $;}\\
4n + 3, & \; \; \; \mbox{when $v \in \{a^i : 1 \leq i \leq 4n - 1, i \ne 2n\} \cup\{a^{2i + 1}b : 1 \leq i \leq 2n - 1  \} $.}
\end{array}\right.\]By the proof of Theorem \ref{detour-ecentricity}, we have the following:
\begin{itemize}
\item $dds_D(e) = (1, 2n, \underbrace{0, 0, \ldots, 0}_{4n - 3}, 1, 0, 6n- 2)$,

\item $dds_D(a^{2n}) = (1, \underbrace{0, 0, \ldots, 0}_{4n -2},1,  2n, 6n -2)$,

\item $dds_D(a^{i}) = (1, \underbrace{0, 0, \ldots, 0}_{4n}, 4n-1, 2n, 2n)$,

\item $dds_D(a^{2i}b) = (1,1, 2n -1, \underbrace{0, 0, \ldots, 0}_{4n -3},1, 0,6n -2)$,

\item $dds_D(a^{2i + 1}b) = (1, \underbrace{0, 0, \ldots, 0}_{4n}, 3,2n, 6n -4)$.
\end{itemize}

Thus, the result holds.
\end{proof}

\begin{proposition}
In $\c$, we have
\begin{enumerate}[\rm (i)]
\item $Int(\mathcal{P}_e(SD_{8n}) = K_2$,

\item $Cl(\mathcal{P}_e(SD_{8n}) = \mathcal{P}_e(SD_{8n})$.
\end{enumerate}
\end{proposition}

\begin{proof}
(i) In view of Lemma \ref{nbd}, $x$ is not a complete vertex  if and only if $x \in \{e,  a^{2n}\}$. By Theorems \ref{SD_8n-boundary-complete} and \ref{SD_8n-boundary-interior}, $x$ is an interior vertex if and only if $x \in \{e, a^{2n}\}$. Thus, we have $Int(\mathcal{P}_e(SD_{8n}) = K_2$.\\
\noindent(ii) In view of Lemma \ref{nbd}, we have 
\[{\rm deg}(v) =\left\{ \begin{array}{ll}
 1, & \; \; \; \mbox{when $v = a^{2i}b$ for some $i$, where $1 \leq i \leq 2n$;}\\
3, & \; \; \; \mbox{when $v = a^{2i + 1}b$ for some $i$, where $0 \leq i \leq 2n - 1$;}\\
4n -1, & \; \; \; \mbox{when $v \in \langle a \rangle \setminus \{e,  a^{2n}\}$;}\\
\end{array}\right.\]

Now, we observe that for any pair of non-adjacent vertices $x$ and $y$,  ${\rm deg}(x) + {\rm deg}(y) \leq 6n + 1 < 8n$. Thus, we have $Cl(\mathcal{P}_e(SD_{8n})) = \mathcal{P}_e(SD_{8n})$.
\end{proof}

\begin{theorem}
In $\c$, we have 
\[Ecc(\mathcal{P}_e(SD_{8n})) = \mathcal{P}_e(SD_{8n}) \setminus \{e\}. \]
\end{theorem}

\begin{proof}
Since $e$ is adjacent to all the vertices of $\mathcal{P}_e(SD_{8n})$ we get $ecc(e) = 1$. For $v \ne e$, we have $d(v, e) = ecc(e)$. Therefore, $v$ is an eccentric vertex in $\mathcal{P}_e(G)$. If $e$ is an eccentric vertex in $\mathcal{P}_e(SD_{8n})$, then $ 1 = d(v, e) = ecc(v)$ for some non-identity vertex says $v$. This implies that $v$ is adjacent to all the vertices of $\mathcal{P}_e(SD_{8n})$. It follows that $v = e$ (see Lemma \ref{nbd}); a contradiction.
\end{proof}

\subsection{Resolving Polynomial} In this subsection, we obtain the resolving polynomial of $\c$. First, we recall some of the basic definitions and necessary results from \cite{kumar2021commuting}.
For $z$ in $\Gamma$, we say that $z$ \emph{resolves} $u$ and $v$ if $d(z, u) \ne d(z, v)$. A subset $U$ of $V(\Gamma)$ is a \emph{resolving set} of $\Gamma$ if every pair of vertices of $\Gamma$ is resolved by some vertex of $U$. The least cardinality of a resolving set of $\Gamma$ is called the \emph{metric dimension} of $\Gamma$ and is denoted by $\operatorname{dim}(\Gamma)$. An $i$-\emph{subset} of $V(\Gamma)$ is a subset of $V(\Gamma)$ of cardinality $i$. Let $\mathcal R(\Gamma, i)$ be the family of resolving sets which are $i$-subsets and $r_i = |\mathcal R(\Gamma, i)|$. Then we define the \emph{resolving polynomial} of a graph $\Gamma$ of order $n$, denoted by $\beta(\Gamma, x)$ as $\beta(\Gamma, x) = \mathop{\sum}_{i= dim(\Gamma)}^{n} r_ix^i$. The sequence $(r_{dim(\Gamma)}, r_{dim(\Gamma) + 1}, \ldots, r_n)$ of coefficients of $\beta(\Gamma, x)$ is called the \emph{resolving sequence}.  Two distinct vertices $u$ and $v$ are said to be \emph{true twins} if N$[u] = $N$[v]$. Two distinct vertices $u$ and $v$ are said to be \emph{false twins} if N$(u) = $N$(v)$. The vertices $u$ and $v$ are said to be \emph{twins} if $u$ and $v$ are either true twins or false twins. For more details on twin vertices, we refer the reader to \cite{a.carmen-metric}. A set $U \subseteq V(\Gamma)$ is said to be a \emph{twin-set} in $\Gamma$ if $u, v$ are twins for every pair of distinct pair of vertices $u, v \in U$. In order to obtain the resolving polynomial $\beta(\c, x)$, the following results will be useful. 

\begin{remark}[{\cite[Remark 3.3]{a.Ali-2016}}]\label{r.twin-set}
If $U$ is a twin-set in a connected graph $\Gamma$ of order $n$ with $|U| = l \geq 2$, then every resolving set for $\Gamma$ contains at least $l-1$ vertices of $U$.
\end{remark}

\begin{proposition}[{\cite[Proposition 3.5]{a.Ali-2016}}]\label{p.resolving-SD_8n}
Let $\Gamma$ be a connected graph of order $n$. Then the only resolving set of cardinality $n$ is the set $V(\Gamma)$ and a resolving set of cardinality $n-1$ can be chosen $n$ possible different ways.
\end{proposition}

%\begin{lemma}\label{nbd-even}
%In $\c$, for $n \in \mathbb{N}$, we have
%\begin{enumerate}[\rm(i)]
%\item  {\rm N}$[x] = \{e,x \}$ if and only if $x \in \{a^{2i} b\}$, where $1 \leq i \leq 2n$.
%\item  {\rm N}$[x] = \langle a \rangle$ if and only if $x \in \langle a \rangle \setminus \{e, a^{2n}\}$.
%\item  {\rm N}$[x] = \{e, a^{2n},x,a^{2n}x \}$ if and only if $x \in H_i=\{a^{2i+1}b,a^{2n+2i+1}b\}$, where $0\leq i\leq n-1.$
%\end{enumerate}
%\end{lemma}

\begin{proposition}\label{p.metric-SD_8n}
The metric dimension of $\c$ is $7n - 4$.
\end{proposition}

\begin{proof}
In view of Lemma \ref{nbd}, we get twin-sets $ \langle a \rangle \setminus \{e,  a^{2n}\}, \{a^{2i}b\}$ where $1 \leq i \leq 2n$ and $\{a^{2j+1}b, a^{2n+2j+1}b\}$ where $0 \leq j \leq n-1$.  By Remark \ref{r.twin-set}, any resolving set in $\c$ contains at least $7n-4$ vertices. Now we provide a resolving set of size $7n-4$. By Lemma \ref{nbd}, notice that the set $R = \{a^{2i}b : 1\leq i \leq 2n-1 \} \cup \{a^i : i \ne 1, 2n,4n \}\cup \{a^{2i+1}b : 0\leq i \leq n-1\}$ is a resolving set of size $7n - 4$. Consequently, ${\rm dim}(\c) = 7n - 4$. %We may now suppose that $n$ is odd. By Lemma \ref{nbd-odd}, note that $\langle a \rangle \setminus \{e, a^n, a^{2n}, a^{3n}\}, \{e, a^n a^{2n}, a^{3n}\}$ and $\{a^ib, a^{n +i}b, a^{2n + i}b, a^{n +3i}b\}$, where $1 \leq i \leq n$, are twin sets in $\c$. In view of Remark \ref{r.twin-set}, any resolving set in $\c$ contains at least $7n-2$ vertices. Further, it is routine to verify that the set $R_{\rm odd} = \{a^ib, a^{n + i}b, a^{2n + i}b : 1 \leq i \leq n\} \cup \{a^i : i \ne 1, 2n \}$ is a resolving set of size $7n-2$. Thus, ${\rm dim}(\c) = 7n - 2$. 
\end{proof}

\begin{theorem}
The resolving polynomial of $\c$ is given below:\\
\[ \beta(\c, x) = x^{8n} + 8n x^{8n-1} +\mathop{\sum}_{i= 7n-4}^{8n - 2} r_ix^i,\]
where \[ r_i =  \left\{
\begin{array}{ll}
     \vspace{0.25cm}
     \sum _{j=0}^{i-(7n-4)}\binom{n}{j}2^{n-j}k_{7n+j-i} & 7n-4\leq i\leq 7n-1 \\
     \vspace{0.25cm}
      \sum _{j=0}^{4}\binom{n}{i-7n+j}2^{n-(i-7n+j)}k_j & 7n\leq i\leq 8n-4 \\
    
       \sum _{j=0}^{8n-i}\binom{n}{n+j-(8n-i)}2^{n-(n+j-(8n-i))}k_j & 8n-3\leq i \leq 8n-2 
\end{array} 
\right. \]

and $k_0=1, k_1=6n, k_2= 8n^2+8n-3, k_3=16n^2-2n-2 ,k_4=8n^2-4n$.
\end{theorem}
\begin{proof}
In view of Proposition \ref{p.metric-SD_8n}, we have dim$(\c) = 7n - 4$.  It is sufficient to find the resolving sequence $(r_{7n - 4}, r_{7n -3}, \ldots, r_{8n -2},  r_{8n -1}, r_{8n})$.  By the proof of Proposition \ref{p.metric-SD_8n}, any resolving set $R$ satisfies the  following:
\begin{itemize}
\item $|R \cap (\langle a \rangle \setminus \{e,  a^{2n}\})| \geq 4n - 3$;
\item $|R \cap  H_l| \geq 1$, where $H_l=\{a^{2l+1}b,a^{2n+2l+1}\}$ for $0 \leq l \leq n-1$;
\item $|R \cap \{a^{2i}b : 1 \leq i \leq 2n\}| \geq 2n-1$. 
\end{itemize}
Let $T= SD_{8n}\setminus \displaystyle \bigcup\limits_{l = 0}^{n-1} H_l $. For $|R| = 7n -4$, choose exactly one element from each $H_l$, $4n-3$ elements from $(\langle a \rangle \setminus \{e,  a^{2n}\})$ and $2n-1$ elements from $\{a^{2i}b : 1 \leq i \leq 2n\}$. Therefore, we have 
$$r_{7n-4}=\binom{2}{1}^n\binom{4n-2}{4n-3}\binom{2n}{2n-1}=2^{n+1}n(4n-2)$$
\\
For $7n-3 \leq |R| = i\leq 7n-1$, with above restriction we have one of the following cases:-
\begin{itemize}
    \item Choose $6n-(7n-i)$ elements from $T$ and $n$ elements from $\bigcup\limits_{l = 0}^{n-1} H_l$.
    
\item Choose $6n-(7n-i)-1$ elements from $T$ and $n+1$ elements from $\bigcup\limits_{l = 0}^{n-1} H_l$. 
\[ \begin{array}{ll} 
	 & \vdots\\
\vspace{-0.75cm}\\	
 & \vdots
\end{array}\]
\item Choose $6n-4$ elements from $T$ and $n+(i-7n+4)$ elements from $\bigcup\limits_{l = 0}^{n-1} H_l$.

\end{itemize}

Thus,
$$r_{i}=\sum _{j=0}^{i-(7n-4)}\binom{n}{j}2^{n-j}k_{7n+j-i}$$ 
where $k_0=1, k_1=6n, k_2= 8n^2+8n-3, k_3=16n^2-2n-2 ,k_4=8n^2-4n$ and $k_t$   is the number of ways of selecting $6n-t$ elements from $T$ in $R$ for 
where, $0 \leq t \leq 4$.
\\
\\
For $7n\leq |R| = i \leq 8n-4$, with above restriction we have one of the following cases:-
\begin{enumerate}[(i)]
\item Choose $6n$ elements from $T$ and $n+(i-7n)$ elements from $\bigcup\limits_{l = 0}^{n-1} H_l$.
\item Choose $6n-1$ elements from $T$ and $n+1+(i-7n)$ elements from $\bigcup\limits_{l = 0}^{n-1} H_l$.
\item Choose $6n-2$ elements from $T$ and $n+2+(i-7n)$ elements from $\bigcup\limits_{l = 0}^{n-1} H_l$.
\item Choose $6n-3$ elements from $T$ and $n+3+(i-7n)$ elements from $\bigcup\limits_{l = 0}^{n-1} H_l$.
\item Choose $6n-4$ elements from $T$ and $n+4+(i-7n)$ elements from $\bigcup\limits_{l = 0}^{n-1} H_l$.
\end{enumerate}
Therefore, we have 
$$r_{i}=\sum _{j=0}^{4}\binom{n}{i-7n+j}2^{n-(i-7n+j)}k_j. $$
 \\

For $  8n-3\leq  |R|=i \leq 8n-2$, with above restriction we have one of the following cases:-
\begin{itemize}
   \item Choose $6n$ elements from $T$ and $n+n-(8n-i)$ elements from $\bigcup\limits_{l = 0}^{n-1} H_l$.
\item Choose $6n-1$ elements from $T$ and $n+n-(8n-i)+1$ elements from $\bigcup\limits_{l = 0}^{n-1} H_l$.
\[ \begin{array}{ll} 
	& \vdots\\
\vspace{-0.75cm}\\	
	& \vdots
\end{array}\]
\item Choose $6n-(8n-i)$ elements from $T$ and $2n$ elements from $\bigcup\limits_{l = 0}^{n-1} H_l$.
\end{itemize}

Thus, we get
$$r_{i}= \sum _{j=0}^{8n-i}\binom{n}{n+j-(8n-i)}2^{n-(n+j-(8n-i))}k_j. $$
%For $|R| = i\geq 7n -4$, there exist $v_1, v_2, \ldots, v_{8n-i} \in SD_{8n} \setminus R$. Therefore we have one of the following:
%\begin{enumerate}[(i)]
%\item $v_j \in \langle a \rangle \setminus \{e, a^{2n}\}$ for some $j$  and \\ $ v_1, v_2, \ldots, v_{j-1}, v_{j+1}, v_{j+2}, \ldots v_{8n-i} \in \left(\displaystyle \bigcup\limits_{i = 1}^{2n} \{a^ib, a^{2n +i}b\}\right)\cup \{e, a^{2n}\}$.
%\item $v_1, v_2, \ldots, v_{8n-i} \in \left(\displaystyle \bigcup\limits_{i = 1}^{2n} \{a^ib, a^{2n +i}b\}\right)\cup \{e, a^{2n}\}$.
%\end{enumerate}  
\\
%For $i = 6n -2$, (ii) does not hold so $v_j \in \langle a \rangle \setminus \{e, a^{2n}\}$ and \\ $ v_1, v_2 \ldots, v_{j-1}, v_{j+1}, v_{j+2}, \ldots, v_{8n-i} \in \left(\displaystyle \bigcup\limits_{i = 1}^{2n} \{a^ib, a^{2n +i}b\}\right)\cup \{e, a^{2n}\}$. Therefore, we obtain $r_{6n-2} = 2^{2n + 1}(4n - 2)$. Now for fixed $i, \; 6n - 1 \leq i \leq 8n - 2$, we get $r_i  = 2^{8n - i} \left\{ \binom{2n + 1}{8n - i} + (2n - 1) \binom{2n + 1}{8n - i -1}  \right\}$.  
By Proposition \ref{p.resolving-SD_8n}, $r_{8n-1}= 8n$ and $r_{8n} = 1$. Hence the result holds.  
\end{proof}

 \section{Laplacian spectrum}
Hamzeh $et \ al.$ \cite{Filoment2017spectrum}, obtained the Laplacian spectrum of the power graph of $SD_{8n}$, when $n=2^{\alpha}$. In this section, we investigate the Laplacian spectrum  of the enhanced power graphs of the semidihedral, generalized quaternion and dihedral groups, respectively. Consequently, we provide the number of spanning trees of $\c,  \d $ and $\a$. 

\begin{theorem}\label{lap-even-SD_8n}  
The characteristic polynomial of the Laplacian matrix of $\c$ is 
\[\Phi(L(\c), x) =x(x -8n)(x-6n)(x -4n)^{4n-3}(x -2)^{n}(x-4)^{n}(x-1)^{2n}.\] 
\end{theorem}

\begin{proof}
The Laplacian matrix $L(\c)$ is the $8n \times 8n$ matrix given below, where the rows and columns are indexed in order by the vertices $e = a^{4n}, a^{2n}, a, a^2, \ldots, a^{2n - 1},a^{2n + 1}, a^{2n + 2}, \ldots, a^{4n - 1}$ and then $ab, a^{2n+1}b,a^3b,a^{2n+3}b,  \ldots,a^{2n-1}b,  a^{4n-1}b,a^2b,a^4b, \ldots, a^4nb=b$.
	\[L(\c)  = \displaystyle \begin{pmatrix}
	8n-1 & -1 & -1 & -1& \cdots \cdots & -1  &-1 & \cdots \cdots &-1  \\
	-1 &  &  & & &  &  &  &  \\
	-1  & &   &  &  &  &  &  &    \\
	-1  &  &   & A &  &    &  &  \mathcal O  & & \\
	\; \; \vdots & \; \;& & &   &  &  &  &  & \\ 
	\; \; \vdots& \; \; & & & &  & &  &  & \\ 
	-1  &  &   &  &  &  & &  & &     \\ 
	-1  &  &   &  &  &  & &  & &     \\
	
	\; \; \vdots& \; \;& & \mathcal O'& &  & & B &  & \\ 
	\; \; \vdots& \; \; & & & &  & &  &  & \\ 
	-1  & &   & &   &  & & & &  \\
	\end{pmatrix}\]
	where
	\[A =\displaystyle \begin{pmatrix}
	6n-1 & -1 & -1 & -1& \cdots \cdots & -1  &-1 & \cdots \cdots &-1  \\
	-1 &  &  & & &  &  &  &  \\
	-1  & &   &  &  &  &  &  &    \\
	-1  &  &   & A_1 &  &    &  &  \mathcal O_1 & & \\
	\; \; \vdots & \; \;& & &   &  &  &  &  & \\ 
	\; \; \vdots& \; \; & & & &  & &  &  & \\ 
	-1  &  &   &  &  &  & &  & &     \\ 
	-1  &  &   &  &  &  & &  & &     \\
	
	\; \; \vdots& \; \;& & \mathcal O_1'& &  & & A_2 &  & \\ 
	\; \; \vdots& \; \; & & & &  & &  &  & \\ 
	-1  & &   & &   &  & & & &  \\
	\end{pmatrix}\] 
	 such that $A_1= 4nI_{4n-2}-J_{4n-2}$ and  \[ A_2= \displaystyle \begin{pmatrix}
		C & O_2 & \cdots \cdots   &O_2 \\
    O_2 & C &  \cdots \cdots   & O_2  \\
    \; \; \vdots & \; \; \vdots & & \; \; \vdots \\ 
	\; \; \vdots&\; \; \vdots & & \; \; \vdots \\ 

	  O_2 & O_2 &  \cdots \cdots   & C      \\

	\end{pmatrix}\] with \[ C= \displaystyle \begin{pmatrix}
	3 &-1 \\
   -1 & 3 
   \end{pmatrix}, \]
	
	 $O_2$ is the zero matrix of order $2\times 2$, $B = I_{2n}$, $\mathcal O$ and  $\mathcal O_1$ are the zero matrices of size $(6n - 1) \times (2n)$ and $(4n - 2) \times (2n)$, respectively, $\mathcal O'$ and $\mathcal O'_1$  are the transpose matrix of $\mathcal O$ and $\mathcal O_1$, respectively. Then the characteristic polynomial of $L(\c)$ is 
	
	\[\Phi(L(\c), x)  = \displaystyle \begin{vmatrix}
	x - (8n-1) & 1 & 1 & 1& \cdots \cdots & 1  & 1 & \cdots \cdots &1  \\
	1 &  &  & &  &   &  &  &  \\
	1  & &   &  &  &  &  &  &    \\
	1  &  &   & (xI_{6n-1} - A) &  &    &  &  \mathcal O  & & \\
	\vdots & & & &   &  &  &  &  & \\ 
	\vdots& 	 & & & &  & &  &  & \\ 
	1  &  &   &  &  &  & &  & &     \\ 
	1  &  &   &  &  &  & &  & &     \\
	
	\vdots&   & &\mathcal O' &  & & & (xI_{2n}-B)  & \\ 
	\vdots&  & & & &  & &  &  & \\ 
	1  & &   & &   &  & & & &  \\
	\end{vmatrix}.\]
	
	Apply row operation $R_1 \rightarrow (x -1)R_1 - R_2 - \cdots -R_{8n}$ and then expand by using first row, we get\\
	
	$\Phi(L(\c), x) = \frac{x(x -8n)}{(x - 1)}\displaystyle \begin{vmatrix}

	    (xI_{6n-1} - A)&   &  &  & \mathcal O    \\

	    \mathcal O'& &  & &(xI_{2n}-B)   \\ 
	   
	\end{vmatrix}$\\
	
%	Again, apply row operation $R_1 \rightarrow (x -2)R_1 - R_2 - R_3 -\cdots - R_{8n -1}$ and then expand by using first row, we get
%	\[\Phi(L(\c), x) = \frac{x(x -8n)^2}{(x-2)}  \displaystyle \begin{vmatrix}
%	xI_{4n -2} - A & \mathcal O\\
%	\mathcal O' & xI_{4n} - B
%	\end{vmatrix}.\] 
By using Schur's decomposition theorem \cite{b.spectra}, we have \[\Phi(L(\c), x) = \frac{x(x -8n)}{(x-1)} |xI_{6n-1} - A|\cdot |xI_{2n}-B|.\] Clearly, 
	$$|xI_{2n}-B| =(x -1)^{2n}.$$ 
	Now we obtain
	\[|xI_{6n-1} - A| =\displaystyle \begin{vmatrix}
	x - (6n-1) & 1 & 1 & 1& \cdots \cdots & 1  & 1 & \cdots \cdots &1  \\
	1 &  &  & & &  &  &  &  \\
	1  & &   &  &  &  &  &  &    \\
	1  &  &   & (xI_{4n-2}-A_1) &  &    &  &  \mathcal O_1 & & \\
	\; \; \vdots & \; \;& & &   &  &  &  &  & \\ 
	\; \; \vdots& \; \; & & & &  & &  &  & \\ 
	1  &  &   &  &  &  & &  & &     \\ 
	1  &  &   &  &  &  & &  & &     \\
	
	\; \; \vdots& \; \;& & \mathcal O_1'& &  & & (xI_{2n}-A_2) &  & \\ 
	\; \; \vdots& \; \; & & & &  & &  &  & \\ 
	1  & &   & &   &  & & & &  \\
	\end{vmatrix}\] \

 Apply row operation $R_1 \rightarrow (x -2)R_1 - R_2 - R_3 -\cdots - R_{6n -1}$ and then expand by using first row, we get
	\[ |xI_{6n-1} - A|= \frac{(x-1)(x -6n)}{(x-2)}  |xI_{4n-2}-A_1|\cdot |xI_{2n}-A_2|.\] 
	
	Clearly, $|xI_{2n}-A_2|=(x-2)^{n}(x-4)^{n}$. Now we obtain $|xI_{4n-2}-A_1|=|xI_{4n-2}-(4nI_{4n-2}-J_{4n-2})|.$ It is easy to compute the characteristic polynomial of the matrix $J_{4n -2}$ is $x^{4n -3}(x-4n + 2)$. It is well known that if $f(x) = 0$ is any polynomial and $\lambda$ is an eigenvalue of the matrix $P$, then $f(\lambda)$ is an eigenvalue of the matrix $f(P)$. Consequently, the eigenvalues of the matrix $A$ are $4n$ and $2$. Note that if $x$ is an eigenvector of $J_n$ corresponding to the eigenvalue $0$, then $x$ is also an eigenvector of the matrix $A$ corresponding to eigenvalue $4n$. Since the dimension of the null space of $J_{4n -2}$ is  $4n - 3$  so that the multiplicity of the eigenvalue $4n$ in the characteristic polynomial of the matrix $A$ is $4n - 3$. Thus, $|xI_{4n -2} - A| = (x-4n)^{4n -3} (x-2)$ and hence the result holds.   
\end{proof}

\begin{corollary}
The Laplacian spectrum of $\c$ is given by 
\[\displaystyle \begin{pmatrix}
0 & 1 & 2 &  4 & 4n &6n& 8n\\
 1 & 2n & n & n & 4n -3 &1 & 1\\
\end{pmatrix}.\]
\end{corollary}

By {\cite[Corollary 4.2]{a.Mohar}}, we have the following corollary. 

\begin{corollary}
The number of spanning trees of $\c$ is $2^{11n - 5}3n^{4n - 2}$.
\end{corollary}

%\begin{remark}
%In \cite{Filoment2017spectrum}, the Laplacian spectrum of the power graph of $SD_{8n}$  for  $n=2^{\alpha}$ was investigated.
%\end{remark}
%$\Phi(L(\c), x) =x(x -8n)^2 (x -4)^{2n}(x -2)^{2n}(x-4n)^{4n -3}.$ 
Hamzeh $et \ al.$ \cite{Filoment2017spectrum}, obtained the Laplacian spectrum of the power graph of $Q_{4n}$ when $n=2^{\alpha}$. Now we determine the Laplacian spectrum of generalized quaternion group $Q_{4n}$ and dihedral group $D_{2n}$, for arbitrary $n\in \mathbb{N}$.
 
\begin{theorem}\label{lap-odd-SD_8n}
The characteristic polynomial of the Laplacian matrix of $\d$ is  \[\Phi(L(\d), x) =x(x -4n)^2(x -4)^{n}(x -2)^{n}(x-2n)^{2n -3}.\]
\end{theorem}
\begin{proof}
The Laplacian matrix $L(\d)$ is the $4n \times 4n$ matrix given below, where the rows and columns are indexed by the vertices $e = a^{2n}, a^{n}, a, a^2, \ldots, a^{n - 1},a^{n + 1},  a^{n + 2}, \ldots, a^{2n - 1},$ and then  $ab, a^{n+1}b,a^2b,a^{n+2}b  \ldots, a^nb, a^{2n}b=b$.
	\[L(\d)  = \displaystyle \begin{pmatrix}
	4n-1 & -1 &  \cdots \cdots & -1  &-1 & \cdots \cdots &-1  \\
	-1 & 4n - 1 &  \cdots \cdots & -1  &-1  & \cdots \cdots &-1  \\
    -1  & -1 &    &    &  &   & & \\
	-1  & -1 &   A &    &  &  \mathcal O  & & \\
	\; \; \vdots & \; \;\vdots&    &  &  &  &  & \\ 
	-1  & -1 &  &  & &  & &     \\ 
	\; \; \vdots& \; \; \vdots& \mathcal O'& &  & B &  &  \\ 
	-1  & -1&    &  & & & &  \\
	\end{pmatrix}\]
	where $A = 2nI_{(2n - 2)} - J_{(2n - 2)}$ and \[B= \displaystyle \begin{pmatrix}
	C & O_2 & \cdots \cdots   &O_2 \\
    O_2 & C &  \cdots \cdots   & O_2  \\
    \; \; \vdots & \; \; \vdots & & \; \; \vdots \\ 
	\; \; \vdots&\; \; \vdots & & \; \; \vdots \\ 

	  O_2 & O_2 &  \cdots \cdots   & C      \\

	\end{pmatrix}\] with \[ C= \displaystyle \begin{pmatrix}
	3 &-1 \\
   -1 & 3 
   \end{pmatrix},\] $O_2$ is the zero matrix of order $2\times 2$, $\mathcal O$ is the zero matrix of order $(2n-2)\times (2n)$ and $\mathcal O'$ is the transpose of  $\mathcal O$.
 Then the characteristic polynomial of $L(\d)$ is 
	
	\[\Phi(L(\d), x)  = \displaystyle \begin{vmatrix}
	x-(4n-1) & 1 &  \cdots & 1  &1 & \cdots 1  \\
	1 &x-(4n-1) &  \cdots  & 1  &1 & \cdots  1  \\

	1  & 1&  &  &  &  &    \\
	1  & 1 &  xI_{2n-2}- A &    &  &  \mathcal O  & \\
	\vdots & \vdots&  &  &  &  &  \\ 

	1  & 1 &   &  & &  &     \\ 
	1  & 1 &    &  & &  &      \\
	
	\vdots&  \vdots& \mathcal O'& &  & xI_{2n} - B &    \\ 
	
	1  & 1&    &  & & & \\
	\end{vmatrix}\]
	
	Apply the following row operations consecutively
	\begin{itemize}
		\item $R_1 \rightarrow (x -1)R_1 - R_2  - \cdots - R_{4n}$
		\item $R_2 \rightarrow (x -2)R_2 - R_3 - \cdots - R_{4n}$
		
	\end{itemize}	
	
	and then expand, we get
	
	\[\Phi(L(\d), x) = \frac{x(x -4n)^2}{(x-2)}  \displaystyle \begin{vmatrix}
	xI_{2n-2} -A & \mathcal O\\
	\mathcal O' & xI_{2n} - B
	\end{vmatrix} = \frac{x(x -4n)^2}{(x-2)} |xI_{2n-2}- A|\cdot|xI_{2n} - B|.\]
	By the similar argument used in the proof of Theorem \ref{lap-even-SD_8n}, we obtain  $$|xI_{2n-2} -A| = (x-2n)^{2n -3} (x-2),\ \ \ |xI_{2n}-B|=(x-2)^{n}(x-4)^{n}$$ 
\\	

%	apply the following row operations consecutively $R_i \rightarrow (x -5)R_i - R_{i +1}  - \cdots - R_{4n}$ where $1 \leq i \leq n$ and then on solving, we get \[|xI-B| = \frac{(x -4)^n (x-8)^n}{(x-5)^n} \displaystyle \begin{vmatrix}
%	(x-7)I_{n} & I_{n} & I_{n}\\
%	I_{n} & (x-7)I_{n} & I_{n}\\
%	I_{n} & I_{n} & (x-7)I_{n}\\
%	\end{vmatrix}.\]

%$\begin{itemize}
%\item For $1 \leq i \leq n$  $R_i \rightarrow (x -6)R_i - R_{i +1}  - \cdots - R_{3n}$
%\item For $n + 1 \leq i \leq 2n$  $R_i \rightarrow (x -7)R_i - R_{i +1}  - \cdots - R_{3n}$
%\end{itemize}
%and then expand,  we obtain \[|xI-B| = \frac{(x -4)^n (x-8)^{3n}}{(x-7)^n} |(x-7)I_n| = (x -4)^n (x-8)^{3n}.\]
Thus, the result holds.
\end{proof}

\begin{corollary}
The Laplacian spectrum of $\d$ is given by 
\[\displaystyle \begin{pmatrix}
0 & 2 &  4 & 2n & 4n\\
 1 & n & n & 2n -3 & 2\\
\end{pmatrix}.\]
\end{corollary}

By {\cite[Corollary 4.2]{a.Mohar}}, we have the following corollary.

\begin{corollary}
The number of spanning trees of $\d$ is $2^{5n - 1} n^{2n - 2}$.
\end{corollary}

\begin{remark}
In the similar lines of $\c$, we get \[\Phi(L(\a), x)=x(x-1)^n(x-n)^{n-2}(x-2n).\] 
\end{remark}

\begin{corollary}
The Laplacian spectrum of $\a$ is given by 
\[\displaystyle \begin{pmatrix}
0 & 1 &  n & 2n \\
 1 & n & n-2 & 1 \\
\end{pmatrix}\]
\end{corollary}   
By {\cite[Corollary 4.2]{a.Mohar}}, we have the following corollary.

\begin{corollary}
The number of spanning trees of $\a$ is $n^{n-2}$.
\end{corollary} 
\section{Concluding Remarks}
We end this manuscript with the important remarks on distant properties and detour distant properties of the enhanced power graph of the generalized quaternion group $Q_{4n}$ and dihedral group $D_{2m}$. F. Ali $et \ al.$ \cite{ali2020power}, obtained the distant properties and detour distant properties of the power graph $\mathcal{P}(G)$, when $G$ is       $Q_{4n}$ or $D_{2m}$, where $n=2^k$ and $m=p^\alpha$, $k, \alpha \in \mathbb{N}$, $p$ is an odd prime. In this case, by \cite[Theorem 28]{a.Aalipour2017} note that $\mathcal{P}(Q_{4n})= \d$ and $\mathcal{P}(D_{2m})=\mathcal{P}_e(D_{2m})$. Moreover, note that, for arbitrary $n$, 
$$\d =K_2\vee (K_{n-2}\cup nK_2)$$
and, for arbitrary $m$
$$\mathcal{P}_e(D_{2m})=K_1 \vee (K_{m-1}\cup \overline{K}_m).$$
Thus, in view of \cite[Proposition 3.1 and Proposition 3.3]{ali2020power}), the results on distant properties and detour distant properties of $\mathcal{P}(Q_{4n})$ and $\mathcal{P}(D_{2m})$ will also hold for $\mathcal{P}_e(Q_{4n})$ and $\mathcal{P}_e(D_{2m})$ for arbitrary $n,m \in \mathbb{N}.$

\section{Acknowledgement}
The first author gratefully acknowledge for providing financial support to CSIR  (09/719(0110)/2019-EMR-I) government of India. The third author wishes to acknowledge the support of MATRICS Grant  (MTR/2018/000779) funded by SERB, India.

\end{document}